\documentclass{aptpub}

\usepackage[dvipsnames]{color}

\authornames{W. Bischoff, A. Gegg} % insert the authors here for use in running head
\shorttitle{Boundary Crossing Probabilities for $(q,d)$-Slepian-Processes} % insert short title here for use in running head

\newcommand{\R}{\ensuremath{\mathbb{R}}} %reelle Zahlen
\renewcommand{\N}{\ensuremath{\mathbb{N}}} %nat Zahlen: 1,2,3,..
\renewcommand{\P}{\ensuremath{\mathbb{P}}}  %vorher: \P war Zeichen f\"{u}r Zeilenumbruch: {\P}

\newcommand{\Deq}{\ensuremath{\stackrel{\mathcal{D}}{=}}}
\newcommand{\dif}{\mathrm{d}} %Abstand beim Integral

\numberwithin{equation}{section}  % If you number theorems, etc. within sections,
\begin{document}

\title{Boundary Crossing Probabilities for (q,d)-Slepian-Processes} % insert title - use \\ if it requires more than one line.

\authorone[Catholic University Eichst\"{a}tt-Ingolstadt]{Wolfgang Bischoff} % Affiliation is just the name of your university or institution

\addressone{Faculty of Mathematics and Geography,
Catholic University Eichst\"{a}tt-Ingolstadt, \\ D-85071 Eichst\"{a}tt,
Germany} % Your postal address goes here.

\author[Catholic University Eichst\"{a}tt-Ingolstadt]{Andreas Gegg} % Affiliation is just the name of your university or institution

\begin{abstract}
For $0<q< d$ fixed let $W^{[q,d]}=(W^{[q,d]}_t)_{t\in {[q,d]}}$ be a \emph{$(q,d)$-Slepian-process} defined as centered, stationary Gaussian process with continuous sample paths and covariance
 \begin{align*}
C_{W^{[q,d]}}(s,s+t) = (1-\frac{t}{q})^+ \quad, q\leq s\leq s+t\leq d.
 \end{align*}
Note that
 \begin{align*}
\frac{1}{\sqrt{q}}(B_t-B_{t-q})_{t\in [q,d]},
\end{align*}
where $B_t$ is standard Brownian motion, is a $(q,d)$-Slepian-process.
In this paper we prove an analytical formula for the boundary crossing probability $\P\left(W^{[q,d]}_t > g(t) \; \text{for some } t\in[q,d]\right)$, $q< d\leq 2q$, in the case $g$ is a piecewise affine function. This formula can be used as approximation for the boundary crossing probability of an arbitrary boundary by approximating the boundary function by piecewise affine functions.\\
\end{abstract}



\section{Introduction} \label{sec:introduction}
The signal plus noise process
$$
Y_t=s(t)+B_t, ~t\in [0,d],
$$
where $d>0$ is a known constant, $s\in L^2[0,d]$ is a deterministic signal and the noise B is standard Brownian motion, is often used in statistics as model for many different situations. By $L^2[0,d]$ we denote the set of squared-integrable functions with respect to the Lebesgue-measure on $[0,d]$.
In order to monitor the process $(s(t)+B_t)_{t\in[0,d]}$ it is convenient to consider the path of the process during a moving window $[t-q,t]$ of length $q$ moving over time $t\in [q,d],$ where  $0<q<d$ is fixed. Using this information of the windows over time $t\in [q,d]$ a decision on a hypothesis of interest has to be made.  Often it is of interest to test the
null-hypothesis
\begin{align}\label{testhyp}
 H_0: s=c, ~~\mbox{where $c \in\R$ is a known or unknown constant}.
\end{align}
A simple statistic for each window is the difference $s(t)+B_t-s(t-q)-B_{t-q}, t\in [q,d],$ of the observations at the two boundary points of the interval $[t-q,t]$. Under $H_0$ and by normalizing the variance to 1 we get the stochastic process
$$
\frac{1}{\sqrt{q}}(B_t-B_{t-q})_{t\in [q,d]}.
$$
This process is at least useful to detect increasing or decreasing signals $s$. Note that the covariance of the above process is given by
   \begin{align}\label{cov}
C(t,t+u) = (1-\frac{t}{q})^+ \quad, q\leq t\leq t+u\leq d,
 \end{align}
where $a^+=a$ if $0\leq a$ and $=0$ if $a<0$. We call a centered, stationary Gaussian process $W^{[q,d]}=(W^{[q,d]}_t)_{ t\in [q,d]}$ with continuous paths and covariance function given in (\ref{cov}) (q,d)-Slepian process.
Slepian (1961) investigated processes of the form $W^{[1,d]}, 1<d$. Nowadays these processes are called Slepian processes, see Azais and Wschebor (2009). Slepian-processes or modifications of them have been studied by several authors. We refer the interested reader to Mehr and Mc Fadden (1965), Shepp (1966, 1971), Bar-David (1975), Ein-Gal and Bar-David (1975),
Cressie (1980, 1981), Abrahams (1984, 1986), Chu et al. (1995), Nikitin and Orsingher (2004/2006), Fuchang and Li (2007), Gegg (2013), Liu et al. (2014), Bischoff and Gegg (2015) and the references given therein.

Boundary crossing probabilities are of specific interest for stochastic processes. For instance,
\begin{equation}
  \label{eq:7}
\P\left( W^{[q,d]}_t > g(t) \; \text{for some } t\in[q,d]\right)
\end{equation}
where $g:[q,d]\to\R$ is a boundary function, is required to establish one-sided tests of Kolmogorov type for the null-hypothesis (\ref{testhyp}).

Note that for $0<q< d$
\begin{align} \label{rel}
W^{[q,d]}\; \Deq \;\frac{1}{\sqrt{q}}(B_t-B_{t-q})_{t\in [q,d]}\; \Deq \;
(B_\frac{t}{q}-B_{\frac{t}{q}-\frac{q}{q}})_{t\in [q,d]},
\end{align}
where ``$\Deq$'' means ``identical in law''. Hence, by putting $e:=\frac{d}{q}, u:=\frac{t}{q}$ and $h(u):=g(uq), u\in[1,e],$
\begin{align} \label{rel2}
\P\left( W^{[q,d]}_t > g(t)\; \text{for some } t\in[q,d]\right)=
\P\left( W^{[1,e]}_u > h(u) \; \text{for some } u\in[1,e]\right).
\end{align}
Therefore the boundary crossing probability for $(q,d)$-Slepian processes can be restricted to boundary crossing probabilities of Slepian processes $W := W^{[1,e]}$. Nevertheless, we state our results for $(q,d)$-Slepian processes. For, the $(q,1)$-Slepian processes appear in models with compact experimental region as described at the very beginning of this introduction and $(2,d)$-Slepian bridges are connected with Brownian bridges, see below.

In the literature analytic formulas of boundary crossing probabilities for Slepian processes are only known if the boundary $g$ is a constant or an affine function, see  Ein-Gal and Bar-David (1975), Abrahams (1984). It seems to be impossible to show such an analytic result for an arbitrary boundary function $g$. Therefore, we suggest the following procedure. At first we present an analytic formula for the boundary crossing probability of a $(q,d)$-Slepian process by considering the boundary function piecewisely. By this formula the boundary crossing probability can be approximated by approximating the boundary function piecewisely by simple functions as, for instance, constant or affine functions. It is important to emphasize that we can show our results for $d\leq 2q$ only, because $(q,d)$-Slepian processes have not the same properties for $d\leq 2q$ and $2q<d$, respectively. All proofs are given in an appendix.
%%%%%%%%%%%%%%%%%%%%%%%%%%%%%%%%%%%%%%%%%%%%%%%%%%%%%%%%%%%%%%%%%%%%%%%%%%%
%%%%%%%%%%%%%%%%%%%%%%%%%%%%%%%%%%%%%%%%%%%%%%%%%%%%%%%%%%%%%%%%%%%%%%%%%%%%%%
\section{Boundary Crossing Probability}

Given an arbitrary boundary function $g$ an approximation of the boundary crossing probability of a $[q,d]$-Slepian-process $W^{[q,d]}$ can be determined by approximating the boundary function $g$ piecewise by simple functions. Wang and P\"otzelberger (1997) applied such an approach to Brownian motion. They used affine functions as simple functions. It is more complicated to apply this idea to Slepian processes because Brownian motion is markovian, Slepian processes are not.

At first we prove an expression for the boundary crossing probability by considering the boundary function piecewisly without simplifying it.
%%%%%%%%%%%%%%%%%%%%%%%%%%%%%%%%%%%%%%%%%%%%%%%%%%%%%%%%%%%%%%%%%%%%%%%
\begin{thm}  \label{th:7}
Let $q,d\in \R$ with $0<q<d\leq 2q$. Let $n+1, n\in\N,$ fixed points in time $t_0,\dots,t_n$ with
$
q=t_0 < t_1 < t_2 < \ldots < t_n=d,
$
be given and let $c:=\frac{q^{(n+1)/2}}{2^{n}\pi^{(n+1)/2}\sqrt{(3q-d)(d-q)}}$. Then for a measurable function $g:[q,d]\to\R$ holds true
\begin{align*}
& \P\left( W_t^{[q,d]} > g(t) \quad \text{for some } t\in[q,d] \right) = 1 - c\cdot\left(\prod_{i=1}^{n-1}\frac{\sqrt{t_{i+1}-q}}{\sqrt{(t_{i+1}-t_{i})(t_i-q)}}\right)\\
&\times\int_{-\infty}^{g(t_0)} \int_{-\infty}^{g(t_n)} \int_{-\infty}^{g(t_{n-1}) } \:
    \cdots \: \int_{-\infty}^{g(t_1)} \;
\exp\left[-\frac{q}{4}\left(\frac{(x_0+x_n)^2}{3q-d}+\frac{(x_0-x_n)^2}{d-q}\right) \right]\\
&\times\P\left(  W^{[q,d]}_t \leq g(t) \; \forall t\in(t_{0},t_{1}) \; | \; W^{[q,d]}_{t_{0}}=x_{0},  W^{[q,d]}_{t_{1}}=x_{1}\right)\\
&\times \prod_{i=1}^{n-1} \P\left(  W^{[q,d]}_t \leq g(t) \; \forall t\in(t_{i},t_{i+1}) \; | \; W^{[q,d]}_{t_{i}}=x_{i},  W^{[q,d]}_{t_{i+1}}=x_{i+1}\right)\\
&\times
 \exp\left[-\frac{q}{4}\left( \frac{ (x_{i}-x_{0})^2}{t_{i}-q}+\frac{ (x_{i+1}-x_{i})^2}{t_{i+1}-t_{i}}-\frac{(x_{i+1}-x_0)^2}{t_{i+1}-q}
  \right)\right] \;\dif x_1 \: \cdots \: \dif x_{n-1}  \dif x_n \dif x_0 \; .
\end{align*}
\end{thm}
%%%%%%%%%%%%%%%%%%%%%%%%%%%%%%%%%%%%%%%%%%%%%%%%%%%%%%%%%%%%%%%%
The proof is given in the appendix \ref{sec:proofs-general}.

For a Slepian process $W^{[1,d]}, 1<d\leq 2,$ Ein-Gal and Bar-David (1975) state the probability
$$
P\left(W^{[1,d]}_t \leq  b_i \; \text{for all } t\in(t_{i-1},t_i), i=1,\ldots,n \; | \;  W^{[1,d]}_1=x_0, W^{[1,d]}_d=x_n\right),
$$
where $x_0, x_n\in\R$, for piecewise constant boundaries $b_i\in \R, i=1,\dots,n,$ as iterated integrals. Our formula given in Theorem \ref{th:7} is an extension of their result, see also the proof of Theorem \ref{th:7}.

We denote
\begin{align} \label{bridge}
\P\left( W^{[q,d]}_t \leq g(t) \; \text{ for all } t\in(t_i,t_{i+1})
\;|  W^{[q,d]}_{t_i}=x_i, W^{[q,d]}_{t_{i+1}}=x_{i+1} \right)
\end{align}
the non-crossing probability of a
\emph{$(q,d)$-Slepian-bridge} starting in $x_i$ at $t=t_i$ and terminating in
$x_{i+1}$ at $t=t_{i+1}$. This probability occurring under the integrals in Theorem \ref{th:7} can be analytically calculated for simple boundary functions $g$.
By the stationarity of $W^{[q,d]}$ we get for $x_1<g(t_1), x_2<g(t_2), q\leq t_1<t_2\leq d$ and $h=t_2-t_1$
\begin{align}\label{hitting}
  & \P\left( W^{[q,d]}_t \leq g(t) \; \text{ for all } t\in(t_1,t_2)\; | \;
  W^{[q,d]}_{t_1}=x_1, W^{[q,d]}_{t_2}=x_2 \right) \nonumber\\
&\quad = \P\left( W^{[q,d]}_{t} \leq g(t_1-q+t) \; \text{ for all } t\in(q,q+h)\; | \;
  W^{[q,d]}_q=x_1, W^{[q,d]}_{q+h}=x_2 \right) \nonumber\\
&\quad = 1-\int_q^{q+h} \pi_{g(t_1-q+\cdot)}\left(t \; | \; W^{[q,d]}_q=x_1, W^{[q,d]}_{q+h}=x_2
\right)\; \dif t  \; ,
\end{align}
where \begin{equation*} \label{eq:10a}
   \pi_{f}^{[q,d]}\left(t \; | \; W^{[q,d]}_q=x_1, W^{[q,d]}_{q+h}=x_{2}\right)
 \end{equation*}
is a Lebesgue-density of the double conditioned first hitting time of $f(s), s\in [q,q+h]$, at $t$ under the condition $W^{[q,d]}_q=x_1, W^{[q,d]}_{q+h}=x_{2}$. Indeed, for  Slepian processes $W^{[1,d]}$ Ein-Gal and Bar-David (1975) determined a Lebesgue-density of the double conditioned first hitting time $\pi_{f}^{[1,d]}\left(t \; | \; W^{[1,d]}_q=x_1, W^{[1,d]}_{q+h}=x_{2}\right)$ for a constant $f=b$. So they obtained the non-crossing probability of a Slepian bridge by the above considerations. Using this formula, see Ein-Gal and Bar-David (1975) page 551, and (\ref{rel}) we get the following result for a ${[q,d]}$-Slepian bridge and a constant boundary.
%%%%%%%%%%%%%%%%%%%%%%%%%%%%%%%%%%%%%%%%%%%%%%%%%%%%%%%%%%%%%%%%%%%%%%%%%%%%%%%%%%%%%%
\begin{lem} \label{lem:3}
Let $q,d\in \R$ with $0<q<d\leq 2q$. For $t_i,t_{i+1}\in[q,d]$ with $t_i<t_{i+1}$,  $b_i\in\R$ and
$x_i,x_{i+1}\in(-\infty,b_i]$ holds
\begin{align*}
\P&\left( W_t^{[q,d]} \leq b_i \; \text{ for all } t\in(t_i,t_{i+1}) \; | \;  W^{[q,d]}_{t_i}=x_i, W^{[q,d]}_{t_{i+1}}=x_{i+1} \right)\\
&\quad = 1- \exp\left( -\frac{q(b_i-x_i)(b_i-x_{i+1})}{t_{i+1}-t_i}  \right)\; .
\end{align*}
\end{lem}
%%%%%%%%%%%%%%%%%%%%%%%%%%%%%%%%%%%%%%%%%%%%%%%%%%%%%%%%%%%%%%%%%%%%%%%%%%%%%%%%%
Note that this probability is closely connected to the
non-crossing probability of a standard Brownian bridge, since for a
standard Brownian motion $B=(B_t)_{t\geq 0}$ holds true for $q\leq t_1< t_1+h\leq d$:
\begin{align*}
\P&\left( W^{[q,d]}_t \leq  b \; \text{ for all } t\in(t_1,t_1+h) \; | \;  W^{[q,d]}_{t_1}= x_1, W^{[q,d]}_{t_1+h}=  x_2 \right)\\
&= 1- \exp\left( -\frac{2( b-  x_1)( b- x_2)}{2h/q}   \right)\\
&=\P\left( B_t \leq  b \; \text{ for all } t\in(t_1,t_1+2h/q) \; | \;  B_{t_1}= x_1, B_{t_1+2h/q}=  x_2 \right).
\end{align*}
For the last equation see e.g. Siegmund (1986) p.~375. This formula shows that the boundary crossing probability of a constant for a $(2,d)$-Slepian bridge and a Brownian bridge coincide if $h\leq 2$.

Finally, we give an analytical formula for (\ref{bridge}) in case g is an affine function. The proof is given in appendix \ref{app 2}.
%%%%%%%%%%%%%%%%%%%%%%%%%%%%%%%%%%%%%%%%%%%%%%%%%%%%%%%%%%%%%%%%%%%%%%%%%%%%%%%%%%%%%%%%
\begin{thm}  \label{th:18}
Let $0<q<d\leq 2q$ and let $a,b\in \R$ be fixed. For $[t_i,t_{i+1}]\subseteq[q,d], h:=t_{i+1}-t_i$ and $x_i<b$, $x_{i+1}<b+ah$ holds true:
\begin{align*}
&\P\left( W_t^{[q,d]} \leq b+a(t-t_i) \; \text{ for all } t\in[t_i,t_{i+1}]\; | \;
  W^{[q,d]}_{t_i}=x_i, W^{[q,d]}_{t_{i+1}}=x_{i+1} \right)\\
 & = 1- \frac{\sqrt{qh}(b-x_i)}{2\sqrt{\pi}}\exp\left(\frac{q(x_{i+1}-x_i)^2}{4h}\right)\times\\
  &\quad\quad\quad\quad\quad \times  \int_0^h \frac{1}{\sqrt{(h-s)s^3}}\exp\left[-\frac{q}{4}\left( \frac{ (b+as-x_{i})^2}{s}+\frac{ (x_{i+1}-b-as)^2}{h-s}
  \right)\right] \; \dif s .
\end{align*}

\end{thm}
%%%%%%%%%%%%%%%%%%%%%%%%%%%%%%%%%%%%%%%%%%%%%%%%%%%%%%%%%%%%%
%%%%%%%%%%%%%%%%%%%%%%%%%%%%%%%%%%%%%%%%%%%%%%%%%%%%%%%%%%%%%
\appendix
 \section{Proofs} \label{secApp:proofs}
For both proofs we need the finite dimensional distributions of the Slepian pocess.
Let $W=W^{[1,e]}$ be a Slepian process, where $e\in\R$ is a constant with $1<e$. Let $\varphi(W_{s_1}=\cdot,\dots,W_{s_m}=\cdot)$ denote a density of the finite dimensional distribution of $(W_{s_1},\dots,W_{s_m})$ and let $\varphi(W_{s_1}=\cdot,\dots,W_{s_m}=\cdot ~|W_{u_1}=a_1,\dots,W_{u_{\ell}}=a_{\ell})$ denote a density of the conditional finite dimensional distribution of $(W_{s_1},\dots,W_{s_m})|W_{u_1}=a_1,\dots,W_{u_\ell}=a_{\ell}, 1\leq s_1<\dots <s_m\leq e, 1\leq u_1<\dots <u_{\ell}\leq e, s_i\neq u_j$. For $m\geq 2$ Slepian (1961) proved the formula
\begin{multline}  \label{finite dim d}
\varphi(W_{s_1}=x_1,\ldots,W_{s_m}=x_m) =  \frac{1}{2^{m-1}\sqrt{\pi^m (2-s_m+s_1)}}
\; \left(\prod_{i=2}^m \frac{1}{\sqrt{s_i-s_{i-1}}} \right)\\
 \times \exp\left[-\frac{1}{4}\left(\frac{(x_1+x_m)^2}{2-s_m+s_1} +    \sum_{i=2}^m\frac{ (x_i-x_{i-1})^2}{s_i-s_{i-1}}
  \right)\right]\;.
\end{multline}
\noindent
Hence, for $1=u_0<u_i<u_{i+1}, u_n=e, x_0,x_i,x_{i+1},x_n\in\R, i=1,\dots,n-1$, we get the following expressions after some calculations
\begin{align}
&\varphi(W_{u_0}=x_{0}, W_{u_i}=x_{i})= \frac{1}{2\pi\sqrt{ (3-u_i)(u_i-1)}}
\exp\left[-\frac{1}{4}\left(\frac{(x_0+x_i)^2}{3-u_i}+\frac{(x_0-x_i)^2}{u_i-1}\right) \right],\label{2dd}\\
&\varphi(W_{u_i}=x_i|W_{u_0}=x_0,W_{u_{i+1}}=x_{i+1}) =
\frac{\varphi\left(W_{1}=x_0, W_{u_i}=x_i,W_{u_{i+1}}=x_{i+1} \right)}{
\varphi(W_{1}=x_0,W_{u_{i+1}}=x_{i+1})} = \nonumber\\
&\frac{\sqrt{u_{i+1}-1}}{2\sqrt{\pi (u_{i+1}-u_{i})(u_i-1)}}
\exp\left[-\frac{1}{4}\left( \frac{ (x_{i}-x_{0})^2}{u_{i}-1}+\frac{ (x_{i+1}-x_{i})^2}{u_{i+1}-u_{i}}-\frac{(x_{i+1}-x_0)^2}{u_{i+1}-1}
  \right)\right]\;.\label{cond fdd}
\end{align}
%%%%%%%%%%%%%%%%%%%%%%%%%%%%%%%%%%%%%%%%%%%%%%%%%%%%%%%%%%%%%%%%%%%%%%%%%%
%%%%%%%%%%%%%%%%%%%%%%%%%%%%%%%%%%%%%%%%%%%%%%%%%%%%%%%%%%%%%
\subsection{Proof of Theorem \ref{th:7}}
\label{sec:proofs-general}
In a first step we apply Fubini's Theorem and get for $u_0=1, u_n=e$
\begin{align*}
  &\P\left( W_u \leq h(u) \quad \text{for all } u\in[1,e] \right)
  \nonumber \\
=&\int_{-\infty}^{h(u_0)} \int_{-\infty}^{h(u_n)}  \; \P\left(
  W_u \leq h(u) \quad \text{for all } u\in(u_0,u_n) \; | \; W_{u_0}=x_0,
  W_{u_n}=x_n\right) \nonumber \\*
&\hspace{5cm}\times \varphi(W_{u_0}=x_0, W_{u_n}=x_n) \; \dif
x_n \dif x_0 \;. % \label{eq:51}
\end{align*}
Next, we use an idea from Ein-Gal and Bar-David (1975). Note that the Slepian process $W$ has the following Markov-like property (see Slepian (1961), Ein-Gal and Bar-David (1975)): let $1\leq s_1 < s_2 \leq e$, then events defined on $[1,s_1)\cup (s_2,e]$ are stochastically independent of events defined on $(s_1,s_2)$ under the condition $W_{s_1}=x_1$ and $W_{s_2}=x_2$. Since this Markov-like property, we have for $x_0\leq h(u_0), x_n \leq h(u_n)$ and $1=u_0<u_{n-1}<u_n=e$:
  \begin{multline}
 \P\left(  W_u \leq h(u) \quad \text{for all } u\in(u_0,u_n) \; | \; W_{u_0}=x_0,
   W_{u_n}=x_n\right) \nonumber\\
 = \; \int_{-\infty}^{h(u_{n-1}) } \P\left(  W_u \leq h(u) \; \forall
  u\in(u_0,u_{n-1}) \; | \; W_{u_0}=x_0, W_{u_{n-1}}=x_{n-1}\right) \\*
\times \P\left( W_{u}\leq h(u) \;  \forall u\in(u_{n-1},u_n) \; | \;
 W_{u_{n-1}}=x_{n-1},  W_{u_n}=x_n\right) \\*
  \times \varphi(W_{u_{n-1}}=x_{n-1}|W_{u_0}=x_0, W_{u_n}=x_n)
\; \dif x_{n-1} \; . %\label{eq:50}
  \end{multline}
Applying this procedure sequentially to the first probability under the integral and using the densities for the finite dimensional (conditional) distributions calculated above leads to the following result.
Let $n+1$ fixed points in time $u_0,\dots,u_n$ with
$
1=u_0 \leq u_1 \leq u_2 \leq \ldots \leq u_n=e,
$
be given and let $c:=\frac{1}{2^{n}\pi^{(n+1)/2}\sqrt{(3-e)(e-1)}}$. Then for $h:[1,e]\to\R$ holds true:
\begin{align*}
  & \P\left( W_u \leq h(u) \quad \text{for all } u\in[1,e] \right) = \\*
&\int_{-\infty}^{h(u_0)} \int_{-\infty}^{h(u_n)} \int_{-\infty}^{h(u_{n-1}) } \:
    \cdots \: \int_{-\infty}^{h(u_1)} \;
c\cdot \exp\left[-\frac{1}{4}\left(\frac{(x_0+x_n)^2}{3-e}+\frac{(x_0-x_n)^2}{e-1}\right) \right]\\*
&\times\P\left(  W_u \leq h(u) \; \forall u\in(u_{0},u_{1}) \; | \; W_{u_{0}}=x_{0}, W_{u_{1}}=x_{1}\right)\\
&\times \prod_{i=1}^{n-1} \P\left(  W_u \leq h(u) \; \forall u\in(u_{i},u_{i+1}) \; | \; W_{u_{i}}=x_{i}, W_{u_{i+1}}=x_{i+1}\right)\frac{\sqrt{u_{i+1}-1}}{\sqrt{(u_{i+1}-u_{i})(u_i-1)}}
 \nonumber \\*
&\times
 \exp\left[-\frac{1}{4}\left( \frac{ (x_{i}-x_{0})^2}{u_{i}-1}+\frac{ (x_{i+1}-x_{i})^2}{u_{i+1}-u_{i}}-\frac{(x_{i+1}-x_0)^2}{u_{i+1}-1}
  \right)\right] \; \dif
x_1 \: \cdots \: \dif x_{n-1}  \dif x_n \dif x_0 \; .
\end{align*}
Hence, the result stated in Theorem \ref{th:7} follows by (\ref{rel}) and (\ref{rel2}).
%%%%%%%%%%%%%%%%%%%%%%%%%%%%%%%%%%%%%%%%%%%%%%%%%%%%%%%%%%%%%%%%%%%%%%%%%%%%%%%
%%%%%%%%%%%%%%%%%%%%%%%%%%%%%%%%%%%%%%%%%%%%%%%%%%%%
\subsection{Proof of Theorem \ref{th:18}}\label{app 2}
Let $W=W^{[1,e]}, 1<e$, be a Slepian process, let $0<\ell\leq e-1$, and let $g(s)=b+a(s-1), s\in [1,1+\ell],$ be an affine boundary function. In the following we need the first hitting time
$\tau_{g}=\inf\{t\geq  \; | \; W_t > g(t) \}$ and a Lebesgue-density of the first hitting time under the condition that the process $W$ starts in $x_1$ at $t=1$. This density is denoted by  $\pi_{g}(\: \cdot \: | W_1=x_1)$. Abrahams (1984), see also Mehr and McFadden (1965), gives an explicit formula of the conditional first hitting time density for an
affine boundary function
  \begin{align}\label{first hitting time density}
   &\pi_{g}\left(t|W_1= x_1\right)\nonumber\\*
   &= \frac{b- x_1}{(t-1) \sqrt{2\pi (t-1) (3-t)}} \cdot
  \exp\left[-\frac{((b+a(t-1))-x_1(2-t))^2}{2(3-t)(t-1)} \right]\nonumber\\*
&= \frac{b-x_1}{t-1} \; \varphi(W_t=b+a(t-1) | W_1=x_1) \; , \quad t\in(1,e],
   \end{align}
where the last equation can be obtained after some calculations using (\ref{2dd}).

By (\ref{hitting}) it is sufficient to determine a Lebesgue-density of the double conditioned first hitting time $\pi_{g}\left(u \; | \; W_1=x_1, W_{1+\ell}=x_2\right)$ at $u$ under the condition $W_1=x_1< g(1)=b, W_{1+\ell}=x_{2}< g(1+\ell)=b+a\ell$.
Ein-Gal and Bar-David (1975) give such a formula if the boundary $g$ is constant. In the following we use their idea and show that it works for an affine boundary function as well.
By Bayes' Theorem
\begin{align*}
&  \pi_{g}(t \; | \; W_1=x_1, W_{1+\ell}=x_{1+\ell}) \\
 &\quad\quad= \; \frac{\varphi\left(W_{1+\ell}=x_{1+\ell}\; | \; W_1=x_1, \tau_{g}=t \right)  \cdot
\pi_{g}(t\; | \;W_1=x_1)}{\varphi(W_{1+\ell}=x_{1+\ell}|W_1=x_1)} \;.
\end{align*}
Note that $1<t<1+\ell$ and
\begin{equation*}
%  \label{eq:177}
\{ \tau_{g}=t\} = \left\{W_s < g(s) \text{ for all } s\in[1,t),
  W_t=g(t)\right\}.
\end{equation*}
Hence, by the Markov-like property of $W$, see the proof of Theorem \ref{th:7}, the expression above coincides with
\begin{equation*}% \label{eq:46}
\frac{\varphi\left(W_{1+\ell}=x_{1+\ell}\; | \; W_1=x_1,  W_t=g(t)  \right)  \cdot
\pi_{g}(t\; | \;W_1=x_1)}{\varphi(W_{1+\ell}=x_{1+\ell}|W_1=x_1)} \;.
  \end{equation*}
Next we apply (\ref{first hitting time density}) and so we obtain  for $1<t=1+u<1+\ell$
\begin{align*}
&  \pi_{1}(1+u \; | \; W_1=x_1, W_{1+\ell}=x_{1+\ell}) = \frac{b-x_1}{u}\varphi\left( W_u=b+au | W_1=x_1, W_{1+\ell}=x_{1+\ell} \right)\\
&=\frac{\sqrt{\ell}(b-x_1)}{2\sqrt{\pi (\ell-u)u^3}}
  \exp\left[-\frac{1}{4}\left( \frac{ (b+au-x_{1})^2}{u}+\frac{ (x_{1+\ell}-b-au)^2}{\ell-u}-\frac{(x_{1+\ell}-x_1)^2}{\ell}
  \right)\right]\;,
\end{align*}
where (\ref{cond fdd}) is used by calculating the last equation. Hence, the result follows by using (\ref{rel2}) and by transforming the density function.
%%%%%%%%%%%%%%%%%%%%%%%%%%%%%%%%%%%%%%%%%%%%%%%%%%%%%%%%%%%%%%%%%%%%%%%%%%%%%%%%%%%%%%%%%%%%%%

\end{document}